\documentclass[12pt,leqno,amsfonts,amscd]{amsart}
\usepackage{epsfig}
\usepackage{amsmath}
\usepackage{amsfonts}
\usepackage{amssymb}
\usepackage{color}

\newtheorem{theorem}{Theorem}[section]

\newtheorem{corollary}[theorem]{Corollary}

\theoremstyle{definition}
\newtheorem{definition}[theorem]{Definition}

\theoremstyle{remark}
\newtheorem{remark}[theorem]{Remark}
\numberwithin{equation}{section}
\newcommand{\calf}{{\mathcal F}}

\newcommand{\RR}{{\mathbb R}}
\newcommand{\CC}{{\mathbb C}}

 \DeclareMathOperator{\psh}{PSH}

\renewcommand{\phi}{\varphi}

\begin{document}
\title[ Finely Plurisubharmonic Functions and Pluripolarity]{Continuity Properties of Finely Plurisubharmonic Functions and
pluripolarity}
\author{Said El Marzguioui}
\address{KdV Institute for Mathematics, Universiteit van Amsterdam, Postbus 94248 1090 GE, Amsterdam, The Netherlands}
\email{s.elmarzguioui@uva.nl}

\author{Jan Wiegerinck}
\address{KdV Institute for Mathematics, Universiteit van Amsterdam, Postbus 94248 1090 GE, Amsterdam, The Netherlands}
\email{j.j.o.o.wiegerinck@uva.nl}
\subjclass[2000]{
%\maketitle \footnote{2000 Mathematics Subject Classification
32U15, 32U05, 30G12, 31C40}
\begin{abstract}
We prove that every bounded finely pluri\-sub\-harmonic function
can be locally (in the pluri-fine topology) written as the
difference of two usual plurisub\-harmonic functions. As a consequence finely pluri\-sub\-harmonic
functions are continuous with respect to the pluri-fine topology.
Moreover we show that $-\infty$ sets of finely plurisubharmonic functions
are pluripolar, hence graphs of finely
holomorphic functions are pluripolar. 
\end{abstract}
\maketitle

%Let $\Omega $ be a bounded open subset of $\mathbb{C}^n$
\section{Introduction}
The fine topology on an open set $\Omega\subset\RR^n$ is the coarsest
topology that makes all subharmonic functions on
$\Omega$ continuous. A finely subharmonic function is defined on a fine domain, it is upper semi-continuous with respect to the
fine topology, and satisfies an appropriate modification of the mean
value inequality. Fuglede \cite{Fu72} proved the following three
properties that firmly connect fine potential theory to classical
potential theory: finely subharmonic functions are finely
continuous (so there is no super-fine topology), all finely polar
sets are in fact ordinary polar sets, and finely subharmonic
functions can be uniformly approximated by subharmonic functions
on suitable compact fine neighborhoods of any point in their
domain of definition. Another continuity result is what Fuglede calls the {\em Brelot Property}, i.e. a finely subharmonic function is continuous on a suitable fine neighborhood of any given point in its domain, \cite[page 284]{Fu88}, see also \cite[Lemma 1]{Fu76}.

Similarly, the pluri-fine topology on $\Omega\subset\CC^n$ is the coarsest topology that makes all plurisubharmonic (PSH) functions on $\Omega$ continuous. In \cite{E-W2} we introduced finely plurisubharmonic functions as plurifinely upper semicontinuous functions, of which the restriction to complex lines is finely subharmonic.
We will prove the analogs of two of the results mentioned
above. Bounded finely plurisubharmonic functions can locally
be written as differences of ordinary PSH functions (cf.~Section 3),
hence finely plurisubharmonic functions are pluri-finely continuous. We also prove a weak form of the Brelot Property.  Next, finely pluripolar sets are shown to be
pluripolar. This answers natural questions posed e.g. by \cite{Mo03}. As a corollary
we obtain that zero sets of finely holomorphic functions of
several complex variables are pluripolar sets. Partial results in
this direction were obtained in \cite{EMW,EdJo06,E-W2}. A final
consequence is Theorem \ref{F-PlfineH} concerning the pluripolar
hull of certain pluripolar sets.

The pluri-fine topology was introduced in \cite{Fu86-1}, and
studied in e.g., \cite{Bed88, BT87, E-W1, E-W2}. In the rest of
the paper we will qualify notions referring to the pluri-fine
topology by the prefix ``$\mathcal{F}$'', to distinguish them from
those pertaining to the Euclidean topology. Thus a compact
$\mathcal{F}$-neighborhood $U$ of $z$ will be a Euclidean compact
set $U$ that is a neighborhood of $z$ in the pluri-fine topology.
%If $z_0$ is a usual interior point of $U$, then there is nothing
%to prove since in this case $f$ is plurisubharmonic at $z_0$. We
%may therefore assume that $z_0$ is not a Euclidean interior point.
\section{Finely plurisubharmonic and holomorphic functions}
There are several ways to generalize the concepts of plurisubharmonic and of holomorphic functions to the setting of the plurifine topology.
See e.g., \cite{Mo03, E-W2}, and in particular \cite{Fu09} where the different concepts are studied and compared. 
\begin{definition}\label{F-PlfineE}
Let $\Omega $ be an $\mathcal{F}$-open subset of $\CC^n$. A
function $f$ on $\Omega$ is called {\em
$\mathcal{F}$-plurisubharmonic} if $f$ is $\mathcal{F}$-upper
semicontinuous on $\Omega$ and if the restriction of $f$ to any
complex line $L$ is finely subharmonic or $\equiv -\infty$ on any
$\mathcal{F}$-connected component of $\Omega\cap L$.

A subset $E$ of $\CC^n$ is called \emph{$\mathcal{F}$-pluripolar}
if for every point $z\in E$ there is an $\mathcal{F}$-open subset
$U\subset \CC^n$ and an $\mathcal{F}$-plurisubharmonic function
($\not \equiv -\infty$) $f$ on $U$ such that $E\cap U
\subset\{f=-\infty\}$.
\end{definition}

Denote by $H(K)$ the uniform closure on $K$ of the algebra of
holomorphic functions in neighborhoods of $K$.
\begin{definition}\label{F-PlfineF} Let $U\subseteq \CC^n$ be $\mathcal{F}$-open. A function $f$ :\
$U$ $\longrightarrow$ $\CC$ is said to be
$\mathcal{F}$-holomorphic if every point of $U$ has a compact
$\mathcal{F}$-neighborhood $K\subseteq U$ such that the
restriction $f|_K$ belongs to $H(K)$.
\end{definition}
\begin{remark}
The functions defined in Definition \ref{F-PlfineE} are called weakly $\calf$-PSH functions in \cite{Fu09}, whereas the functions in Definition \ref{F-PlfineF} are called strongly $\calf$-holomorphic functions. In \cite{Fu09} strongly $\calf$-PSH functions (via approximation) and weakly $\calf$-holomorphic functions (via holomorphy on complex lines) are defined and it is shown
that the strong properties imply the weak ones.
\end{remark}

The original definition of finely subharmonic functions involves
sweeping-out of measures. If one wants to avoid this concept, one can use the next theorem  as an alternative definition.
\begin{theorem}[Fuglede \cite{Fu74,Fu82}] \label{Diff-th5} A function $\varphi$ defined in an $\mathcal{F}$-open set
$U \subseteq \mathbb{C}$ is finely subharmonic if and only if every
point of $U$ has a compact $\mathcal{F}$-neighborhood $K \subset
U$ such that $\varphi|_{K}$ is the uniform limit of usual
subharmonic functions $\varphi_n$ defined in Euclidean
neighborhoods $W_n$ of $K$.
\end{theorem}
%The following result of Bedford and Taylor (see page 5, and
%Theorem 3.4 in \cite{Bed82}) will be used in the proof of Theorem
%\ref{Brelotpsh}.
%\begin{theorem}\label{Monge-Am} Let $u_j \in \psh(\Omega)\cap L^{\infty}_{loc}(\Omega)$ be %a sequence decreasing to $u
%\in \psh(\Omega)\cap L^{\infty}_{loc}(\Omega)$. Then for every
%compact subset $K$ of $\Omega$ and $\delta >0$
%\begin{equation}
%\lim_{j\rightarrow \infty} C(K\cap \{u_{j}>u+\delta\}, \Omega)=0.
%\end{equation}
%\end{theorem}
%Here $L^{\infty}_{loc}(\Omega)$ denotes the set of locally bounded
%function on $\Omega $, and $C$ denotes the Monge-Amp\`{e}re
%capacity.

Recall also the following property, cf.~\cite{BT87}, which will be
used in the proof of Theorem \ref{PlfineC} and its corollary.
\begin{theorem}(Quasi-Lindel\"{o}f property)\label{F-PlfineG} An
arbitrary union of $\mathcal{F}$-open subsets of $\CC^{n}$ differs
from a suitable countable subunion by at most a pluripolar set.
\end{theorem}

\section{Continuity of Finely PSH Functions}
\begin{theorem}\label{PlfineA} Let $f$ be a bounded
$\mathcal{F}$-pluri\-sub\-harmonic function in a bounded
$\mathcal{F}$-open subset $U$ of $\mathbb{C}^n $. Every point $z
\in U$ then has an $\mathcal{F}$-neighborhood $\mathcal{O} \subset
U$ such that $f$ is representable in $\mathcal{O}$ as the
difference between two locally bounded plurisubharmonic functions
defined on some usual neighborhood of $z$. In particular $f$ is $\calf$-continuous.
\end{theorem}
\begin{proof}%[Proof of Theorem \ref{PlfineA}. ]
 We may assume that $-1<f < 0$ and that $U$
is relatively compact in the unit ball $B(0, 1)$. Let $V \subset
U$ be a compact $\mathcal{F}$-neighborhood of $z_0$. Since the
complement $\complement V$ of $V$ is pluri-thin at $z_0$, there
exist $ 0<r<1$ and a plurisubharmonic function $\varphi$ on
$B(z_0,r)$ such that
\begin{equation}\label{Plfine1}
\limsup_{z \to z_{0},z\in \complement V }\varphi(z)<
\varphi(z_{0}).
\end{equation}
Without loss of generality we may suppose that $\varphi $ is
negative in $B(z_{0}, r)$ and
\begin{equation}\label{Plfine2}
\varphi(z)=-1 \  \text{if} \ z\in B(z_0,r) \backslash V \
\text{and} \ \varphi(z_0)=-\frac{1}{2}.
\end{equation}
Hence
\begin{equation}\label{Plfine4}
f(z)+\lambda \varphi(z) \leq -\lambda  \ \text{for}\ \text{any} \
z \in U\cap B(z_0, r)\backslash V \ \text{and} \ \lambda >0.
\end{equation}
Now define a function $u_{\lambda}$ on $B(z_0, r)$ as follows

\begin{equation}\label{SaidFuG} u_{\lambda}(z)=
    \begin{cases}
        \max\{-\lambda , \ f(z)+\lambda \varphi(z)\} & \text{if $z\in U \cap B(z_0, r)$,}\cr
        -\lambda &\text{if $z\in B(z_0, r)\backslash V$.}
    \end{cases}
\end{equation}
This definition makes sense because $[U \cap B(z_0, r)]\bigcup
[B(z_0, r)\backslash V ]=B(z_0, r)$, and the two definitions of $u_\lambda$ agree
on $U\cap B(z_{0}, r)\backslash V$ in view of (\ref{Plfine4}).

Clearly, $u_{\lambda}$ is $\mathcal{F}$-plurisubharmonic in $U
\cap B(z_0, r)$ and in $B(z_0, r)\backslash V$, hence in all
$B(z_0, r)$ in view of the sheaf property, cf.~\cite{E-W2}. Since
$u_{\lambda}$ is bounded in $B(z_0, r)$, it follows from
\cite[Theorem 9.8]{Fu72} that $u_{\lambda}$ is subharmonic on each
complex line where it is defined. It is a well known result that a
bounded function which is subharmonic on each complex line where
it is defined, is plurisubharmonic, cf.~\cite{Lelong45}. In other
words $u_{\lambda}$ is plurisubharmonic in $B(z_0, r)$.

Since $ \varphi(z_0)=-\frac{1}{2}$, the set
$\mathcal{O}=\{\varphi>-3/4\}$ is an
$\mathcal{F}$-neighborhood of $z_0$, and because $\varphi= -1$ on
$B(z_0, r)\backslash V$, it is clear that $\mathcal{O}\subset V
\subset U$.

Observe now that $-4 \leq \ f(z)+ 4\varphi(z)$, for every $z\in
\mathcal{O}$. Hence
\begin{equation}
f(z)=u_{4}(z)-4\varphi(z), \ \text{for} \ \text{every} \ z\in
\mathcal{O}. 
\end{equation}
We have shown that $f$ is $\calf$-continuous on a neighborhood of each point in its domain, hence $f$ is $\calf$-continuous.
\end{proof}
The proof is inspired by \cite[page 88-90]{Fu72}. 

\begin{corollary}\label{PlfineB} Every $\mathcal{F}$-plurisubharmonic function is
$\mathcal{F}$-continuous.
\end{corollary}
\begin{proof}
%[Proof of Corollary \ref{PlfineB}. ] 
Let $f$ be $\mathcal{F}$-plurisubharmonic
in an $\mathcal{F}$-open subset $\Omega$ of $\mathbb{C}^n$. Let $d<c\in \RR$. The set $\Omega_c=\{f<c\}$ is $\calf$-open. The function $\max\{f, d\}$ is bounded $\calf$-PSH on $\Omega_c$, hence $\calf$-continuous. 
Therefore the set $\{d<f<c\}$ is $\calf$-open, and we conclude that $f$ is $\calf$-continuous.
\end{proof}
The following result gives a partial analog to the Brelot
property. We recall the definition of the {\em relative extremal function} or{\em pluriharmonic measure} of a subset $E$ of an open set $\Omega$, cp. \cite{Bed82,K91}
 %\cite[page 284]{Fu88}, see also \cite[Lemma 1]{Fu76}.

\begin{equation}\label{Plfine3}
U=U_{E,\Omega}=\sup\{\psi \in \text{PSH}^{-}\Omega: \ \psi \leq -1 \
\text{on} E\}.
\end{equation}
It is well known that the upper semi-continuous regularization of
$U$, i.e. $U^{*}(z)=\limsup_{\Omega\ni v \rightarrow z
}U(v)$ is plurisubharmonic in $\Omega$.
\begin{theorem}\label{Quasi-Brelot}(Quasi-Brelot property) Let $f$ be a plurisubharmonic function in the unit ball
$B\subset \CC^n$. Then there exists a pluripolar set $E\subset B$
such that for every $z \in B \setminus E$ we can find an
$\mathcal{F}$-neighborhood $\mathcal{O}_z \subset B$ of $z$ such
that $f$ is continuous in the usual sense in $\mathcal{O}_z$
\end{theorem}
\begin{proof} Without loss of generality we may assume that $f$ is continuous near
the boundary of $B$. By the quasi-continuity theorem (cf.~\cite[Theorem 3.5.5]{K91} and the remark that follows it, see also \cite{Bed82}) we can
select a sequence of relatively compact open subset $\omega_n$ of
$B$ such that the Monge-Amp\`{e}re capacity $C(\omega_n,
B)<\frac{1}{n}$, and $f$ is continuous on $B \setminus \omega_n$.
Denote by $\tilde \omega_n$ the $\mathcal{F}$-closure of
$\omega_n$.

The pluriharmonic measure $U^{*}_{\omega_n,
B}$ is equal to the pluriharmonic measure $U^{*}_{\tilde \omega_n,
B}$, because for a PSH function $\varphi$ the set $\{\varphi \leq
-1\}$ is $\mathcal{F}$-closed, thus $\varphi|_{\omega_n} \leq -1
\Rightarrow \varphi|_{\tilde \omega_n} \leq -1$. Now, using
\cite[Proposition 4.7.2]{K91}
\begin{equation}\label{e1}
C(\omega_n, B)=C^*(\omega_n, B)=\int_\Omega (dd^cU^{*}_{\omega_n,
B})^n=\int_\Omega (dd^cU^{*}_{\tilde \omega_n,
B})^n=C^*(\tilde \omega_n, B).
\end{equation}
Let $E=\bigcap_{n}\tilde \omega_n$. By \eqref{e1}, $C^*(E,
B)\leq C^*(\tilde \omega_n, B)\leq \frac{1}{n}$, for every $n$.
Hence $E$ is a pluripolar subset of $B$.

Let $z \not \in E$. Then there exists $N$ such that $z \not \in
\tilde \omega_N$. Clearly, the set $B\setminus \tilde \omega_N$ is
an $\mathcal{F}$-neighborhood of $z$. Since $f$ is continuous on
$B\setminus \omega_N$, it is also continuous on the smaller set
$B\setminus \tilde \omega_N$ ($\subset B\setminus \omega_N$).
\end{proof}
\begin{remark}
The above Quasi-Brelot property holds also for
$\mathcal{F}$-pluri\-sub\-harmonic functions, in view of Theorem
\ref{PlfineA}.
\end{remark}

\section{$\mathcal{F}$-Pluripolar Sets and Pluripolar Hulls}
In this section we prove that $\mathcal{F}$-pluripolar sets are
pluripolar and apply this to pluripolar hulls.
\begin{theorem}\label{PlfineC} Let $f$ :\ $U$ $\longrightarrow$
$[-\infty,+\infty[$ be an $\mathcal{F}$-plurisubharmonic function
$(\not \equiv -\infty)$ on an $\mathcal{F}$-open and
$\mathcal{F}$-connected subset $U$ of $\CC^{n}$. Then the set
$\{z\in U: \ f(z)=-\infty \}$ is a pluripolar subset of $\CC^{n}$
\end{theorem}

\begin{proof}[Proof of Theorem \ref{PlfineC}.] We may assume that $f<0$. Let $z_{0} \in U$, which we can assume relatively compact in $B(0,1)$.
 We begin by showing that $z_0$ admits an $\mathcal{F}$-neighborhood $W_{z_0}$  such that
$\{f=-\infty\}\cap W_{z_0}$ is pluripolar. If $z_0$ is a Euclidean interior point of $U$, then $f$ is PSH on a neighborhood of $z_0$ and there is nothing to prove. 

If not we proceed as in the proof of Theorem
\ref{PlfineA}. 
Thus, let $V\subset U$ be a compact
$\mathcal{F}$-neighborhood of $z_0$, and $\varphi$ a negative PSH function on $B(z_0, r)$ such that
\begin{equation}\label{Plfine2a}
\varphi(z)=-1 \  \text{if} \ z\in B(z_0,r) \backslash V \
\text{and} \ \varphi(z_0)=-\frac{1}{2}.
\end{equation}

Let $\Phi=U_{B(z_0,r)\setminus V, B(z_0,r)}$ be the pluriharmonic measure defined in \eqref{Plfine3}. By
\eqref{Plfine2a}, we get $\varphi \leq \Phi \leq \Phi^{*}$. In
particular $ -\frac{1}{2}\leq \Phi^{*}(z_0)$.

 Let $f_{n}=\frac{1}{n}\max(f, -n)$. Then
$-1\leq f_{n}<0$. We  define functions $v_{n}(z)$ on $B(z_0, r)$ as
follows.
\begin{equation}\label{Plfine6}
v_{ n}(z)=
    \begin{cases}
        \max\{-1, \ \frac{1}{4}f_{n}(z)+ \Phi^{*}(z)\} & \text{if $z\in U \cap B(z_0, r)$,}\cr
        -1 &\text{if $z\in B(z_0, r)\backslash V$.}
    \end{cases}
\end{equation}
Since $v_n$ is analogous to the function $u_{\lambda}$ in
(\ref{SaidFuG}), the argument in the proof of Theorem \ref{PlfineA} shows that $v_n \in
\psh(B(z_0, r))$. Now for
$z\in U$ such that $f(z)\neq-\infty$ the sequence $f_{n}(z)$ increases to $0$. Thus
$\{v_{n}\}$ is an increasing sequence of PSH-functions. Let $\lim v_n=\psi$. The upper semi-continuous regularization $\psi^{*}$ of
$\psi$ is plurisubharmonic in $B(z_0, r)$. It is a result of \cite{Bed82}, see also Theorem
4.6.3 in \cite{K91}, that the set $E=\{\psi
\neq \psi^{*}\}$ is a pluripolar subset of $B(z_0, r)$.

We claim that $\psi^{*}=\Phi^{*}$ on $B(z_0, r)$. Indeed, $\psi\le\psi^*\le\Phi^*$ because the $v_n$ belong to the defining family \eqref{Plfine3} for $\Phi$.
Now observe that $\psi = \Phi^{*}$ on $B(z_0, r) \setminus
\{f=-\infty\}$, because $v_n=\Phi^{*}= -1$ on $B(z_0, r)\backslash
V$. Hence
\begin{equation}\label{fInTe}
\{ \psi^{*}\ne \Phi^{*} \} \subset \ B(z_0, r) \cap \{f=-\infty\}.
\end{equation}
Clearly, the set $\{ \psi^{*} \neq
\Phi^{*}\}$ is  $\mathcal{F}$-open. In view of Theorem 5.2 in
\cite{E-W2} it must be empty because it is contained in the $-\infty$-set of a finely plurisubharmonic function.

Let $z\in \{\Phi^{*}>-\frac{2}{3} \}\cap\{f=-\infty\}$. Then it follows 
from the definition of $v_{n}$ and the claim that 
$$\psi (z) = -\frac{1}{4} +
\Phi^{*}(z) = -\frac{1}{4} + \psi^{*}(z).$$
Thus $z\in E$. 
%is shows that
%$\{\Phi^{*}>-\frac{2}{3} \} \cap \{f=-\infty\}$ is a subset of $E$. 
Now $\{\Phi^{*}>-\frac{2}{3} \}$ is an
$\mathcal{F}$-neighborhood of $z_{0}$. The conclusion is that every point
$z\in U$ has an $\mathcal{F}$-neighborhood
$W_{z} \subset U$ such that $W_{z}\cap
\{f=-\infty\}$ is a pluripolar set. ( If $f(z)\neq -\infty$ we
could have chosen $W_{z}$ such that $W_{z}\cap
\{f=-\infty\}=\emptyset$.) 

By the Quasi-Lindel\"{o}f property, cf.~Theorem \ref{F-PlfineG} there is a sequence $\{z_n\}_{n\geq1}
\subset U$ and a pluripolar subset $P$ of $U$ such that
\begin{equation}\label{FpLiNd1}
U= \cup_{n}\mathcal{O}_{z_n}\cup P.
\end{equation}
Hence
\begin{equation}\label{FpLindP}
\{f=-\infty\} \subset ( \cup_n\mathcal{O}_{z_{n}} \cap
\{f=-\infty\})\cup P.
\end{equation}
This completes the proof since a countable union of pluripolar
sets is pluripolar.
\end{proof}
\begin{remark} Corollary \ref{PlfineB} and Theorem \ref{PlfineC} give affirmative
answers to two questions in \cite{Mo03}.
\end{remark}
A weaker formulation of Theorem \ref{PlfineC}, but perhaps more
useful, is as follows.
\begin{corollary} Let $f$ :\ $U$ $\longrightarrow$
$[-\infty,+\infty[$ be a function defined in an
$\mathcal{F}$-domain $U \subset \CC^n$. Suppose that every point
$z\in U$ has a compact $\mathcal{F}$-neighborhood $K_{z} \subset
U$ such that $f|_{K_{z}}$ is the decreasing limit of usual
plurisubharmonic functions in Euclidean neighborhoods of $K_z$.
Then either $f\equiv -\infty $ or the set $\{f= -\infty\}$ is
pluripolar subset of $U$.
\end{corollary}
As a byproduct we get the following corollary which recovers and
generalizes the main result in \cite{EMW} to functions of several
variables.
\begin{corollary}\label{PlfineD} Let $h$ :\ $U$ $\longrightarrow$
$\CC$ be an $\mathcal{F}$-holomorphic function on an
$\mathcal{F}$-open subset $U$ of $\CC^n$. Then the zero set of $h$
is pluripolar. In particular, the graph of $h$ is also pluripolar.
\end{corollary}
\begin{proof}[Proof of Corollary \ref{PlfineD}. ] Let $a \in U$. Definition \ref{F-PlfineF} gives us
a compact $\mathcal{F}$-neighborhood $K$
of $a$ in $U$, and a sequences $(h_{n})_{n \geq 0 }$, of
holomorphic functions defined in Euclidean neighborhoods of $K$
such that
$$
h_{n}|_{K} \longrightarrow h|_{K}, \ \text{uniformly}.
$$
For $k\in \mathbb{N}$ we define $v_{n, k}=\max(\log|h_{n}|, -k)$
and $v_{k}=\max(\log|h|, -k)$. Clearly, $v_{n, k}$ converges
uniformly on $K$ to $v_{k}$ as $n\to\infty$. Accordingly,
$v_{k}$ is $\mathcal{F}$-plurisubharmonic on the
$\mathcal{F}$-interior $K'$ of $K$. Since $v_{k}$ is decreasing,
the limit function $\log|h|$ is $\mathcal{F}$-plurisubharmonic in
$K'$. Theorem \ref{PlfineC} shows that the set $K'\cap \{h=0\}$ is
pluripolar. The corollary follows now by application of the
Quasi-Lindel\"{o}f property.
\end{proof}

The pluripolar hull $E_{\Omega}^*$ of a pluripolar set $E$
relative to an open set $\Omega $ is defined as follows.
$$
E_{\Omega}^{*}=\bigcap \{z\in \Omega \ : u(z)= - \infty\},
$$
where the intersection is taken over all plurisubharmonic
functions defined in $\Omega$ which are equal to $-\infty$ on $E$.

The next theorem improves on Theorem 6.4 in
\cite{E-W2}.
\begin{theorem}\label{F-PlfineH} Let $U \subset \CC^{n}$ be an $\mathcal{F}$-domain, and
let $h$ be $\mathcal{F}$-holomorphic in $U$. Denote by
$\Gamma_{h}(U)$ the graph of $h$ over $U$, and let $E$ be a
non-pluripolar subset of $U$. Then $ \Gamma_{h}(U)\subset
(\Gamma_{h}(E))^{\ast}_{\CC^{n+1}}$.
\end{theorem}

\begin{proof} By Corollary \ref{PlfineD} the set $\Gamma_{h}(E)$ is pluripolar subset of $\CC^{n+1}$.
Let $\varphi$ be a plurisubharmonic function in $\CC^{n+1}$ with $
\varphi \not \equiv -\infty $ and $\varphi(z, h(z))=-\infty $, for
every $z \in E$. The same arguments as in the proof of Lemma 3.1
in \cite{EMW} show that the function $z\mapsto \varphi(z, h(z))$
is $\mathcal{F}$-plurisubharmonic in $U$. Since $E$ is not
pluripolar, it follows from Theorem \ref{PlfineA} that $\varphi(z,
h(z))=-\infty$ everywhere in $U$. Hence $ \Gamma_{h}(U)\subset
(\Gamma_{h}(E))^{\ast}_{\CC^{n+1}}$.
\end{proof}

\section{Some further questions}

\textbf{Question 1} Let $f$ be an $\mathcal{F}$-plurisubharmonic
function defined in an $\mathcal{F}$-open set $U \subseteq \CC^2$.
Suppose that for each point $z\in U$ there is a compact
$\mathcal{F}$-neighbourhood $K_z$ such that $f$ is continuous (in
the usual sense) on $K_z$. Is it true that $f|_{K_z}$ is the
uniform limit of usual plurisubharmonic functions $\varphi_n$
defined in Euclidean neighborhoods $W_n$ of $K_z$?.

\textbf{Question 2} It is also interesting to figure out whether
the assumption in the above question is automatically fulfilled.
This would be the Brelot property for
$\mathcal{F}$-plurisubharmonic function.

Many other questions remain open. For example, we do not know  the
answer to the following.

\textbf{Question 3} Is this concept of an
$\mathcal{F}$-plurisubharmonic function biholomorphically
invariant?

\end{document}